\documentclass[12pt,reqno]{elsart}
\usepackage{amssymb}

\usepackage{amsmath}

\def\qed{{\unskip\nobreak\hfil\penalty50\hskip1em\nobreak\hfil {\vrule
height7ptwidth7ptdepth0pt}
\parfillskip=0pt\finalhyphendemerits=0\par}\vskip3mm}

\newtheorem{theorem}{Theorem}
\newtheorem{lemma}[theorem]{Lemma}

\theoremstyle{definition}

\renewcommand{\marginpar}[1]{}
\newcommand{\al}{\alpha}

\def\obox {\hbox{\vrule\vbox{\hrule
 \hbox spread 4pt{\hfil\vbox spread 6pt{\vfil}\hfil}
\hrule}\vrule}}

\def\qed{{\unskip\nobreak\hfil\penalty50\hskip1em\nobreak\hfil {
\obox}
\parfillskip=0pt\finalhyphendemerits=0\par}\vskip3mm}

\newcommand{\proof}{\emph{Proof }\,}

\def\bR{\hbox{$\mathbb R$}}

\def\L1o{\hbox{$L^1(\omega)$}}
\def\Lio{\hbox{$L^\infty(\omega^{-1})$}}
\def\Lioo{\hbox{$L^\infty(\omega^{-1}\times \omega^{-1})$}}
\def\LUC1oo{\hbox{$LUC_1(\omega^{-1}\times \omega^{-1})$}}
\def\Mo{\hbox{$M(G,\omega)$}}

\def\tensor{\hbox{$\widehat{\otimes}$}}

\begin{document}
\begin{frontmatter}
 \title
{Amenability of Beurling algebras, corrigendum to a result in ``Generalised notions of amenability, II''}

\author[one]{Fereidoun Ghahramani\thanksref{grant1}}
\ead{Fereidoun.Ghahramani@umanitoba.ca}
\author[two]{Richard J. Loy
}
\ead{Rick.Loy@anu.edu.au}
\author[three]{Yong Zhang\corauthref{cor}\thanksref{grant3}}
\corauth[cor]{Corresponding author.}
\ead{Yong.Zhang@umanitoba.ca}

\address[one]{ Department of Mathematics,
University of Manitoba, Winnipeg R3T 2N2, Canada}
\address[two]{ Mathematical Sciences Institute,
Australian National University, ACT 0200, Australia}
\address[three]{ Department of Mathematics,
University of Manitoba, Winnipeg R3T 2N2, Canada}

\thanks[grant1]{Supported by NSERC grant 2017-05476.} 
\thanks[grant3]{Supported by NSERC grant 2016-238949.}



\begin{abstract}
We fix a gap in the proof of a result in our earlier paper ``Generalised notions of amenability, II'' (J. Funct. Anal. 254 (2008) 1776-1810), and so provide a new proof to a characterization of amenability for Beurling algebras. The result answers a question raised by M. C. White (Bull. London. Math. Soc. 23 (1991) 375–380).
\end{abstract}
\end{frontmatter}
\maketitle

Let $G$ be a locally compact group. A weight function $\omega$ on $G$ is a positive continuous real-valued function on $G$ that satisfies $\omega(xy) \leq \omega(x) \omega(y)$ for all $x,y \in G$. Denote the unit element of $G$ by $e$. Under the condition ${\omega(e) = 1}$,  N.~Gr\o nb\ae k showed in \cite{GRO} that the weighted group convolution algebra $L^1(G, \omega)$, called the Beurling algebra of $G$ for $\omega$, is amenable if and only if $G$ is an amenable group and the function $\Omega(g) = \omega(g)\omega(g^{-1})$ is bounded on $G$. M. C. White raised the question in \cite{white} as to whether the condition $\omega(e) =1$ is removable.  We gave a new proof to Gr\o nb\ae k's result in \cite{GLZ} based on the statement that if $L^1(G, \omega)$ is amenable, then there is a virtual diagonal $u\in (L^1(G, \omega)\tensor L^1(G, \omega))^{**}$ such that $\delta_g\cdot u\cdot \delta_{g^{-1}} = u$ for all $g\in G$. We point out here that the latter condition holds if and only if $L^1(G)$ has a quasi-central bounded approximate identity, that is, a bounded approximate identity $(e_\al)\subset L^1(G)$ satisfying $\delta_g *e_\al - e_\al *\delta_g \to 0$, for all $g\in G$. Of course, this is the case when $G$ is an amenable or SIN group, but is not true in general (see \cite{L-R}). 
This oversight compromises propositions 8-12 of our paper \cite{GLZ}, resulting in them being unavailable in general.

In this note we remove the use of the above additional assumption on $u$ and thus fix the gap in our proof of amenability result for Beurling algebras given in \cite{GLZ}. This gives a complete answer to the above mentioned question of White.

Since $L^1(G, \omega)$ is an ideal of $\Mo$, we may regard it as a Banach $\Mo$-bimodule with the convolution $*$ as the left and the right module multiplications. This naturally makes $L^1(G, \omega)\tensor L^1(G, \omega)$ a Banach $\Mo$-bimodule with module actions defined by
\[
\mu\cdot (a\otimes b) = (\mu*a)\otimes b \quad \text{and}\quad (a\otimes b)\cdot \mu = a\otimes (b*\mu),
\] 
for $\mu \in \Mo$ and $a,b\in L^1(G, \omega)$. The induced $\Mo$-module actions on $(L^1(G, \omega)\tensor L^1(G, \omega))^* = L^\infty(G\times G,\omega^{-1}\times \omega^{-1})$ 
are given by the following formulas. For $\mu\in \Mo$ and $F\in L^\infty(G\times G,\omega^{-1}\times \omega^{-1})$,
\[
[\mu\cdot F](x,y) = \int_G F(x, ys)\,d\mu(s) = \int_G F(x, s)\,d(\delta_y*\mu)(s)\quad (x,y\in G)
\]
\[
[F\cdot\mu ](x,y) = \int_G F(sx, y)\,d\mu(s) = \int_G F(s, y)\,d(\mu*\delta_x)(s) \quad (x,y\in G).
\]

In the sequel, to simplify notation we will denote $L^1(G, \omega)$, $L^\infty(G, \omega^{-1})$ and $L^\infty(G\times G, \omega^{-1}\times \omega^{-1})$ by $\L1o$, $\Lio$ and $\Lioo$, respectively.

\begin{lemma}\label{Deltabdd}
Let $\omega$ be a weight function on a locally compact group $G$. Suppose that $\L1o$ is amenable.  Then the function
$\Omega(g) =\omega(g) \omega(g^{-1})$ is bounded on $G$.
\end{lemma}

\proof  
Let $u\in (\L1o\tensor \L1o)^{**}$ be a virtual diagonal for $\L1o$, so that $f\cdot u = u\cdot f$ and $\pi^{**}(u)\cdot f = f$ for all $f\in \L1o$, where $\pi^{**}$ is the second dual of the multiplication operator $\pi$: $\L1o\tensor \L1o \to \L1o$ defined by $\pi(a\otimes b) = a*b$ for $a,b\in \L1o$. 
Let $(e_i)$ be a bounded approximate identity for $\L1o$ such that $e_i(x)\geq 0$ ($x\in G$) and $\|e_i\|\leq C$ for all $i$. For each $g\in G$ we have,
\[
(e_i*\delta_g)\cdot u \cdot (\delta_{g^{-1}}*e_i) = (u\cdot (e_i*\delta_g))\cdot (\delta_{g^{-1}}*e_i) = u\cdot e_i^2.
\]
Going to a subnet if necessary, we may assume wk*-$\lim _i u\cdot e_i^2$ exists. Define $\Tilde u = \text{wk*-}\lim_i u\cdot e_i^2$. We then have $\pi^{**}(\Tilde u) \cdot f =f$ for $f\in \L1o$ and
\[
\text{wk*-}\lim_i (e_i*\delta_g)\cdot u \cdot (\delta_{g^{-1}}*e_i) = \Tilde u
\]
for all $g\in G$.
Now take a $f\in \L1o$ such that $f$ has compact support $K$ and 
$\int_G f(x) dx = 1$.  Since ${\bf 1}_K\in \Lio$, which is a Banach $\L1o$-bimodule, $F = f\cdot {\bf 1}_K \in  \Lio$, where $\cdot$ denotes the $\L1o$-module multiplication. 
So $\pi^*(F) \in \Lioo$. Indeed
\[
\pi^*(F)(x, y) = F(xy) = \int {\bf 1}_K(xy\xi)f(\xi)d\xi\quad (x,y\in G).
\]
It follows that $\pi^*(F)(x, y) = 0$ whenever $xy\not\in KK^{-1}$.  Set $E = KK^{-1}$,  a compact subset of $G$.
Then
 \begin{equation}\label{nonzero}
\langle\pi^*(F),\Tilde u\rangle = \langle F,\pi^{**}(\Tilde u)\rangle =\langle{\bf 1}_K,\pi^{**}(\Tilde u)\cdot f\rangle
=\langle{\bf 1}_K,f\rangle = \int_K f(x)dx =1\,.
\end{equation}

Define
\begin{equation}\notag
A = \{(x, y) : xy \in E\}\,.
\end{equation}
Then $\pi^*(F)$ has its support contained in $A$, so $\pi^*(F) = \pi^*(F){\bf 1}_A$.

Given $\alpha > 0$, define
\begin{align*}
A_\alpha =&\ \{(x, y)\in A : \omega(x)\omega(y) < \alpha\}\,,\\
B_\alpha = &\ A\setminus A_\alpha = \{(x, y)\in A : \omega(x)\omega(y) \geq \alpha\}\,.
\end{align*}
Clearly  $\pi^*(F) = \pi^*(F){\bf
1}_{A_\alpha}+\pi^*(F){\bf 1}_{B_\alpha}$, and both $\pi^*(F){\bf 1}_{A_\alpha}$ and $\pi^*(F){\bf 1}_{B_\alpha}$ belong to $\Lioo$.

Now we estimate
\begin{align*}
|\langle \pi^*(F){\bf 1}_{B_\alpha}, \Tilde u\rangle| \leq &\ \|\ \pi^*(F){\bf 1}_{B_\alpha}\|\cdot \|\Tilde u\|\\
=&\ \|\Tilde u\| \sup_{B_\alpha} \left| \frac{\pi^*(F)(x,y)}{\omega(x)\omega(y)}\right|\\
=&\ \|\Tilde u\| \sup_{B_\alpha} \left| \frac{F(xy)} {\omega(xy)} \cdot\frac{\omega(xy)}{\omega(x)\omega(y)}\right|\\
\leq&\  \alpha^{-1}\|\Tilde u\|\,\|F\| C_1\,,
\end{align*}
where $C_1 = \sup_{t\in E} \omega(t)$.  It follows that 
\begin{equation}\label{Balpha}
\lim_{\alpha \to \infty}\langle \pi^*(F){\bf 1}_{B_\alpha}, \Tilde u\rangle = 0\,.
\end{equation}

Furthermore, for each $g\in G$,
\begin{align}\notag
|\langle \pi^*(F){\bf 1}_{A_\alpha}, \Tilde u\rangle| = &\lim_i\ |\langle \pi^*(F){\bf 1}_{A_\alpha},(e_i*\delta_g)\cdot u\cdot (\delta_{g^{-1}}*e_i\rangle|\\ 
\leq&\ \limsup_i\|u\| \ \|(\delta_{g^{-1}}*e_i)\cdot \pi^*(F){\bf 1}_{A_\alpha}\cdot (e_i*\delta_g)\|. \label{B estimate}
\end{align}
For $(x,y)\in G\times G$,
\begin{align*}
&\left[(\delta_{g^{-1}}*e_i)\cdot \pi^*(F){\bf 1}_{A_\alpha}\cdot (e_i*\delta_g)\right](x,y) \\
=&\ \int_{G\times G}\,[\delta_{yg^{-1}}*e_i](t)\,\left[\pi^*(F){\bf 1}_{A_\alpha}\right](s,t)\, [e_i*\delta_{gx}](s)\, dt\, ds \\
=& \int_{A_{\al}}[\delta_{yg^{-1}}*e_i](t)\,[\pi^*(F)](s,t) \,[e_i*\delta_{gx}](s)\, dt\, ds \\
=& \int_{A_{\al}}[\delta_{yg^{-1}}*e_i](t)\,F(st) \,[e_i*\delta_{gx}](s)\, dt\, ds.
\end{align*}
So that
\begin{align*}
&\left|\left[(\delta_{g^{-1}}*e_i)\cdot \pi^*(F){\bf 1}_{A_\alpha}\cdot (e_i*\delta_g)\right](x,y)\right|\\ 
\leq&\ \|F\|\ \int_{A_{\al}}[\delta_{yg^{-1}}*e_i](t)\,[e_i*\delta_{gx}](s)\,\omega(st) \, dt\, ds\\
\leq&\ C_1\|F\|\ \int_{A_{\al}}[\delta_{yg^{-1}}*e_i](t)\,[e_i*\delta_{gx}](s)\, dt\, ds.
\end{align*}
We now multiply and divide the last term by $\omega(g)\omega(g^{-1})$.  Note that, for $(s,t)\in A_\al$,
\begin{align*}
\omega(g)\omega(g^{-1}) \leq&\ \omega(gy^{-1}t)\,\omega(sx^{-1}g^{-1})\,\omega(t^{-1}y)\omega(xs^{-1})\\
\leq&\ \omega(gy^{-1}t)\,\omega(sx^{-1}g^{-1})\,\omega^2((st)^{-1})\,\omega(s)\omega(t)\,\omega(x)\omega(y)\\
\leq&\ \al \ C_2^2 \ \omega(gy^{-1}t)\ \omega(sx^{-1}g^{-1})\ \omega(x)\omega(y),
\end{align*}
where $C_2 = sup \{\omega(\tau): \tau\in E^{-1}\}$. 
Thus
\begin{align}\notag
&\|(\delta_{g^{-1}}*e_i)\cdot \pi^*(F){\bf 1}_{A_\alpha}\cdot (e_i*\delta_g)\| \\ \notag
=&\ \sup_{x,y\in G}\left|\frac{\left[(\delta_{g^{-1}}*e_i)\cdot \pi^*(F){\bf 1}_{A_\alpha}\cdot (e_i*\delta_g)\right](x,y)}{\omega(x)\omega(y)}\right|\\ \notag
\leq&\ \frac{\al \ C_1 \ C_2^2 \ \|F\|}{\omega(g)\omega(g^{-1})} \int_{A_{\al}}[\delta_{yg^{-1}}*(e_i\omega)](t)\,[(e_i\omega)*\delta_{gx}](s)\, dt\, ds\\ \label{B part}
\leq&\ \frac{\al \ C_1 \ C_2^2 \ \|F\|C^2}{\omega(g)\omega(g^{-1})} .
\end{align}

Suppose the result is false.  Then there is a sequence $(g_n)\subset G$ such that \break
$\lim_{n\to\infty}\omega(g_n)\omega(g_n^{-1}) = \infty$, whence it follows from (\ref{B estimate}) and (\ref{B part}) that for each $\alpha > 0$,
\begin{equation}\label{Aalpha}
|\langle \pi^*(F){\bf 1}_{A_\alpha}, \Tilde u\rangle| = 0\,.
\end{equation}

Putting (\ref{Balpha}) and (\ref{Aalpha}) together, we may conclude
\[
 \langle \pi^*(F),\Tilde u\rangle = 0\,,
\]
which contradicts (\ref{nonzero}). Therefore, the function $\Omega(g)$ must be bounded on $G$.  \qed

Let $LUC(\omega^{-1})$ be the space of all continuous functions $f\in \Lio$ such that the left translation  $s\mapsto l_sf$ is continuous from $G$ into $\Lio$. 
It is readily seen that $LUC(\omega^{-1})$ is a right $\L1o$ submodule of $\Lio$ and $\L1o$ has a bounded right approximate identity for $LUC(\omega^{-1})$. So 
\[
LUC(\omega^{-1}) = \Lio\cdot \L1o
\]
due to the Cohen factorization theorem. Under the condition $\omega(e) = 1$ it is shown  in \cite{GRO} (see also \cite[Proposition~7.17]{D-L}) that 
$$LUC(\omega^{-1}) = \{f\in \Lio: f/\omega \in LUC(G)\}.$$


We define  $\LUC1oo$ to be  the space of all continuous functions $f\in \Lioo$ that are left uniformly continuous with respect to the first variable, i.e.
\begin{align*}
\LUC1oo = \{& f\in \Lioo\cap C(G\times G): \\
& s\mapsto l_{(s,e)}f: G\to \Lioo \text{ is continuous}\}.
\end{align*}
When $\omega =\bf 1 $, we will denote such space by $LUC_1(G\times G)$. 
It is easy to check that  $|f|\in \LUC1oo$ if $f\in \LUC1oo$. So $\max\{f_1, f_2\}$, $\min\{f_1, f_2\} \in \LUC1oo$ whenever $f_1, f_2\in \LUC1oo$ are real-valued. Moreover, $\LUC1oo$ is a Banach right $\L1o$-submodule of $C(\omega^{-1}\times \omega^{-1})$ and $\L1o$ has a bounded right approximate identity for $\LUC1oo$, where $C(\omega^{-1}\times \omega^{-1}) = \Lioo\cap C(G\times G)$ is regarded as a $\L1o$-submodule of $\Lioo$.  Therefore, due to the Cohen factorization theorem, 
 \[
 \LUC1oo = C(\omega^{-1}\times \omega^{-1}) \cdot \L1o.
 \]

We then note that $\pi^*(LUC(\omega^{-1})) \subset \LUC1oo$. In particular, $\pi^*(\omega)\in \LUC1oo$ if $\omega\in LUC(\omega^{-1})$. The latter condition is automatically satisfied if $\omega(e) = 1$.

Recall that $\phi$: $G\to \bR$ is a character if $\phi(xy) = \phi(x)\phi(y)$ ($x,y\in G$). The following result is due to M. C. White \cite[Lemma~1]{white}.
\begin{lemma}\label{white char}
Let $G$ be an amenable group and $\omega$ be a weight on $G$. Then there is a continuous positive character $\phi$ on $G$ such that
\[ 
\phi(g) \leq \omega(g)\quad (g\in G)\,.
\]
\end{lemma}

We now prove a similar result, replacing the amenability condition on $G$ by amenability condition on $\L1o$.

\begin{lemma}\label{char} Let $\omega$ be a weight on $G$ such that $\omega\in LUC(\omega^{-1})$. Suppose that $L^{\,1}(\omega)$ is amenable.  Then there is a continuous positive character $\phi$ on $G$ such that
\[ 
\phi(g) \leq \omega(g)\quad (g\in G)\,.
\]
\end{lemma}

\proof
If $(u_\al)\in \L1o\tensor \L1o$ is a bounded approximate diagonal for $\L1o$, then so is $(Re(u_\al))$. So there is $u\in (L^1_R(\omega)\tensor L^1_R(\omega))^{**}$ such that $\varphi\cdot u = u\cdot \varphi$ and $\pi^{**}(u)\cdot \varphi = \varphi$ for all  $\varphi\in L^1_R(\omega)$, where $ L^1_R(\omega)$ denotes the real Beurling algebra for the weight $\omega$. Due to this fact, in the following discussion we may simply assume all function spaces involved are real-valued function spaces. In particular $\L1o = L^1_R(\omega)$.

Let $\hat u = u|_{\LUC1oo}$. Then $\hat u\in (\LUC1oo)^*$ and 
\begin{equation}\label{u hat}
\delta_{g^{-1}}\cdot \hat u\cdot \delta_g = \hat u \quad (g\in G)
\end{equation}
since $\LUC1oo =C(\omega^{-1}\times \omega^{-1})\cdot \L1o$. 

For $f\in \LUC1oo^+$ we define
\[
\widetilde{u}(f) = \sup\{  \langle \hat u, \psi\rangle :  |\psi | \leq f, \psi \in \LUC1oo\}\,.
\]
Then $\widetilde{u} \not\equiv 0$ on $\LUC1oo^+$. In fact, 
$\pi^*(\omega)\in \LUC1oo^+$ and
\begin{align}\label{c}
 \widetilde{u}(\pi^*(\omega)) &\geq \langle \hat u, \pi^*(\omega)\rangle = \langle \pi^{**}(u), \omega\rangle = \lim_i \langle \pi^{**}(u), \omega\cdot e_i\rangle \notag\\ 
 & =  \lim_i \langle e_i, \omega\rangle \geq \liminf_i \|e_i\|_\omega \geq 1,
\end{align}
where $(e_i) \subset \L1o^+$ is a bounded approximate identity.

It is standard to check the following affine properties of $\widetilde{u}$: 
\[
\widetilde{u}(cf) = c \widetilde{u}(f), \quad \widetilde{u}(f_1 + f_2) = \widetilde{u}(f_1) + \widetilde{u}(f_2)
\] 
for $c\geq 0$ and $f, f_1, f_2 \in \LUC1oo^+$. For example, one can verify $\widetilde{u}(f_1 + f_2) \leq \widetilde{u}(f_1) + \widetilde{u}(f_2)$ as follows: If $\psi \in \LUC1oo$ satisfies $|\psi| \leq f_1 + f_2$, then we let 
\[
\psi_1 = \max\{-f_1, \min\{\psi, f_1\}\} \text{ and } \psi_2 = \max\{\psi - f_1, \min\{\psi + f_1, 0\}\}.
\]
We have $\psi_1, \psi_2 \in \LUC1oo$, $|\psi_1|\leq f_1$, $|\psi_2|\leq f_2$, and $\psi_1 + \psi_2 = \psi$. These lead to 
\[
\langle \hat u, \psi\rangle = \langle \hat u, \psi_1\rangle +  \langle \hat u, \psi_2\rangle\leq \widetilde{u}(f_1) + \widetilde{u}(f_2)
\]
for all $\psi\in \LUC1oo$ with $|\psi|\leq f_1 + f_2$. So $\widetilde{u}(f_1 + f_2) \leq \widetilde{u}(f_1) + \widetilde{u}(f_2)$ holds. The opposite inequality is obvious. Thus the claimed equality holds.

Clearly, $0 \leq\widetilde{u}(f) \leq \|u\|\, \|f\|$ for $f\in \LUC1oo^+$.  Then $\widetilde{u}$ can be extended to a bounded linear functional on $\LUC1oo$, still denoted by $\widetilde u$, in the obvious manner.  We have $\widetilde{u}\not = 0$, $\|\widetilde u\|\leq \|u\|$ and  $\langle \widetilde{u}, f\rangle
\geq 0$ for $f\in \LUC1oo^+$. Therefore, $\widetilde{u}$ is monotonic. Moreover, given $g\in G$ and $f\in \LUC1oo^+$,it is true that $ \psi \in \LUC1oo$ and $|\psi | \leq f$ if and only if $\delta_g\cdot \psi \cdot \delta_{g^{-1}}
\in\LUC1oo$ and $|\delta_g\cdot \psi \cdot \delta_{g^{-1}}| \leq \delta_g\cdot f \cdot \delta_{g^{-1}}$.
This together with (\ref{u hat}) ensures  
\begin{equation}\label{gu=ug}
\delta_{g^{-1}}\cdot \widetilde{u}\cdot \delta_g = \widetilde{u} \quad (g\in G) .
\end{equation}

By Lemma~\ref{Deltabdd}, there is $M > 0$ such that $\omega(g^{-1})\omega(g) \leq M$ for all $g\in G$. Denote by $\Gamma$ the collection of all finite subsets of $G$. For each $F\in \Gamma$, define
\[
\omega_F (x) = \max_{g\in F}\omega(g^{-1}xg) \quad (x\in G).
\]
We have $\omega_F\in LUC(\omega^{-1})$ and 
\begin{equation}\label{M}
\frac{1}{M}\omega \leq \omega_F \leq M\omega.
\end{equation}

For $g\in G$, let 
\[
W_g(x,y) = \frac{\omega(gx)\omega(gy^{-1})}{\omega(x)\omega(y^{-1})} \quad (x,y\in G).
\]
Then $W_g\in LUC_1(G\times G)$ and
\[
0<\frac{1}{\omega(g^{-1})^2}\leq W_g \leq \omega(g)^2.
\]
Therefore, $\log W_g \in LUC_1(G\times G)$ for each $g\in G$. We note that, with the pointwise multiplication, $LUC_1(G\times G)\LUC1oo \subset \LUC1oo$. So
\[
A^{(F)}_g : =\left( \frac{1}{2}\log W_g\right)\pi^*(\omega_F)\in \LUC1oo
\]
for each $g\in G$ and each $F\in \Gamma$. Furthermore
\begin{equation}\label{bound for A}
-\left(\log \omega(g^{-1})\right)\pi^*(\omega_F)\leq A^{(F)}_g \leq \left(\log \omega(g)\right)\pi^*(\omega_F)\,.
\end{equation}
We note that
\[
W_{g_1g_2}(x,y) = W_{g_1}(g_2x, yg_2^{-1})W_{g_2}(x,y)\,.
\]
We have
\begin{align*}
A^{(F)}_{g_1g_2}(x,y)  &=\left(\left(\frac{1}{2}\log W_{g_1}(g_2 x, yg_2^{-1})\right)
+\left(\frac{1}{2}\log W_{g_2}( x, y)\right)\right)\pi^*(\omega_F)(x,y)\\
    & = \left(\delta_{g_2^{-1}}\cdot A^{(g_2 F)}_{g_1} \cdot \delta_{g_2}\right)(x,y) + A^{(F)}_{g_2}(x,y). 
\end{align*}
To get the second equality above we have used  $\delta_h\cdot \pi^*(\omega_F)\cdot\delta_{h^{-1}} =  \pi^*(\delta_h\cdot \omega_F\cdot\delta_{h^{-1}})$ ( $h\in G$) and 
\[
\delta_h\cdot \omega_F\cdot\delta_{h^{-1}}(x) =\omega_F(h^{-1}xh) =  \max_{g\in F}\omega(g^{-1}(h^{-1}xh)g) = \max_{g\in hF}\omega(g^{-1}xg) = \omega_{hF}(x)\,.
\]
Applying (\ref{gu=ug}) we derive
\begin{equation}\label{sum for A}
\langle \widetilde{u},A^{(F)}_{g_1g_2}\rangle = \langle \widetilde{u}, A^{(g_2F)}_{g_1}\rangle + \langle \widetilde{u}, A^{(F)}_{g_2}\rangle.
\end{equation}
 The net $(\pi^*(\omega_F))_{F\in \Gamma}$ is a bounded increasing net in the commutative unital C*-algebra $\LUC1oo$ whose product is given by the formula  
 $$\Phi\cdot \Psi (x,y) = \frac{\Phi(x,y)\Psi(x,y)}{\omega(x)\omega(y)}.$$
  So the net converges to some $\varUpsilon\in \LUC1oo^{**}$ in the weak* topology of $\LUC1oo^{**}$. From (\ref{c}) and (\ref{M}) we have
\[
\langle\widetilde{u}, \varUpsilon\rangle \geq \frac{1}{M} >0.
\]
Now regard $LUC_1(G\times G)$ as a Banach algebra with the pointwise multiplication. Then $\LUC1oo$ is a Banach $LUC_1(G\times G)$-bimodule (also with pointwise multiplication as the module action). The induced left module action on $\LUC1oo^{**}$ is weak* continuous. So we have 
\[
\text{wk*-}\lim_F A^{(F)}_g = 
 \frac{1}{2}\left(\log W_g)\right) \varUpsilon\,.
\]
Denote the right side by $A_g$. Taking limit in (\ref{sum for A}) we then derive
\[
\langle \widetilde u, A_{g_1 g_2}\rangle =\langle \widetilde u, A_{g_1}\rangle + \langle \widetilde u, A_{g_2}\rangle \quad (g_1, g_2\in G)\,.
\]
On the other hand, by (\ref{bound for A})
\[
-\log\omega(g^{-1})\,\langle \widetilde u, \Upsilon\rangle \leq \langle \widetilde u, A_{g}\rangle \leq \log\omega(g)\,\langle \widetilde u, \Upsilon\rangle \,.
\]
To conclude, we define
\[
\phi(g) = \exp\left(\langle\widetilde u, A_{g}\rangle/\langle\widetilde u, \Upsilon\rangle\right)\quad (g\in G)\,.
\]
Then $\phi$ is a character on $G$ and it satisfies
\[
\frac{1}{\omega(g^{-1})} \leq \phi(g) \leq \omega(g) \quad (g\in G)\,.
\]
The inequality together with the continuity of $\omega$ ensures that $\phi$ is locally bounded, which then implies that $\phi$ is continuous (see the proof of \cite[Lemma~1]{white}). The proof is complete.
\qed

\begin{lemma}\label{equivalent}
Let $\omega$ be a  weight on $G$. Then there is a weight $\bf w$ on $G$ such that ${\bf w}\in LUC({\bf w}^{-1})$ and $\bf w$ is eqivalent to $\omega$, that is, there are $m, M >0$ such that
\begin{equation}\label{equiv}
m\omega(x) \leq {\bf w}(x) \leq M \omega(x) \quad (x\in G).
\end{equation}
\end{lemma}

\proof
The proof is essentially the same as that of \cite[Theorem~3.7.5]{RAS}. 
Let $K$ be a symmetric compact neighborhood of $e$ in $G$. Take a function $\varphi\in C(G)^+$ such that supp$(\varphi)\subset K$ and $\int_G \varphi =1$. Define $\rm w= \omega\cdot \varphi$, i.e.
\[
{\rm w}(x) = \int_G \varphi(\xi)\omega(\xi x)d\xi \quad (x\in G).
\]
Then $\rm w \in LUC(\omega^{-1})^+$. Let $N = \max\{\omega(x): x\in K\}$. 
Then ${\bf w}(x) = N^3{\rm w}(x)$ is a weight on $G$, ${\bf w} \in LUC(\omega^{-1})$, and the equivalence relation (\ref{equiv}) holds for $m=N^2$ and $M = N^4$. The latter two conditions, in turn, imply that ${\bf w} \in LUC({\bf w}^{-1})$. 
 The proof is complete.
\qed

We are now ready to prove our final theorem.

\begin{theorem} Let $\omega$ be a  weight on $G$. Then
$\L1o$ is amenable if and only if $G$ is amenable and the function $\Omega(g)=\omega(g)\omega(g^{-1})$  ($g\in G$) is bounded.
\end{theorem}

\proof 
From Lemma~\ref{equivalent} we may assume $\omega\in LUC(\omega^{-1})$.

Suppose that $\L1o$ is amenable. Then $\Omega$ is bounded due to Lemma~\ref{Deltabdd}. We show that $G$ is amenable. Of course this can be simply done by using the isomorphism $\L1o \cong L^1(G)$ as we did in \cite[Proposition~8.10]{GLZ}. But here we give a direct proof  by  showing that there is a left invariant mean on $LUC(G)$. This proof itself should be of independent interest. Let $\phi$ be the character obtained in Lemma~\ref{char}. Then $\phi\in LUC(\omega^{-1})$, and so $f\phi\in LUC(\omega^{-1})$  for each $f\in LUC(G)$. We have $f\phi\times \phi\in \LUC1oo$. Let $\widetilde u$ be the bounded linear functional on $\LUC1oo$ obtained in the proof of Lemma~\ref{char}. Define
\[
m_0(f) = \langle \widetilde u, f\phi\times \phi\rangle \quad (f\in LUC(G)).
\]
Since $\phi$ is a character we have 
\[
l_gf(x) \,\phi(x)\,\phi(y) = [l_g(f\phi)](x)\,[r_{g^{-1}}\phi](y).
\]
So $(l_gf)\phi\times\phi  = \delta_{g^{-1}}\cdot (f\phi\times \phi)\cdot \delta_g$.
Therefore, 
\[
m_0(l_gf) =  \langle \widetilde u,  \delta_{g^{-1}}\cdot (f\phi\times \phi)\cdot \delta_g\rangle = \langle \widetilde u, f\phi\times \phi\rangle = m_0(f),
\]
i.e. $m_0$ is left invariant. Moreover, $|m_0(f)| \leq  \langle \widetilde u, \phi\times \phi\rangle \|f\|_\infty$ for all $f\in LUC(G)$. Thus 
\[
\|m_0\| \leq  \langle \widetilde u, \phi\times \phi\rangle = m_0(\bf 1).
\]
This implies $\|m_0\| = m_0(\bf 1)$. On the other hand, $m_0 \neq 0$. To see this, we note $\omega(x)\omega(x^{-1}) \leq M$ for some $M>0$ and then
\begin{equation}\label{control omega}
\omega(x) = \omega(x)\phi(x^{-1})\phi(x) \leq \omega(x)\omega(x^{-1})\phi(x) \leq M\phi(x) \quad (x\in G).
\end{equation}
So
\[
 m_0({\bf 1}) = \langle \widetilde u, \phi\times \phi\rangle =  \langle \widetilde u, \pi^*(\phi)\rangle =  \langle \pi^{**}(\widetilde u),  \phi\rangle \geq \frac{1}{M}\langle \pi^{**}(\widetilde u),  \omega\rangle \geq \frac{1}{M}.
\]
All the above show that $m = \frac{m_0}{\|m_0\|}$ is a left invariant mean on $LUC(G)$. So $G$ is amenable.

For the converse, if $G$ is amenable, then by Lemma~\ref{white char} there is a positive character $\phi$ on $G$ such that $\phi \leq \omega$. If in addition $\omega(g)\omega(g^{-1})\leq M$ on $G$, then (\ref{control omega}) holds.
Therefore, $\L1o$ is isomorphic to $L^1(G)$ through the Banach algebra isomorphism $f\mapsto f\phi$. Thus $\L1o$ is amenable since $L^1(G)$ is amenable for the amenable group $G$.
\qed

\textbf{Acknowledgement}: The authors would like to thank the referee for carefully reading the paper and valuable comments.

\bibliographystyle{amsplain}

\end{document}